\documentclass[12pt]{amsart}
\usepackage{amsmath,amstext,amsfonts,amsthm,amssymb,hyperref, amscd,wasysym,stmaryrd}
\usepackage{eucal}
\usepackage[usenames]{color}
\usepackage[all]{xy}                                        %
\swapnumbers
\newtheorem{thm}[subsubsection]{Theorem}
\newtheorem{lem}[subsubsection]{Lemma}
\newtheorem{prp}[subsubsection]{Proposition}
\newtheorem{crl}[subsubsection]{Corollary}

{  \theoremstyle{definition}

}

\newcommand{\Cat}{\mathtt{Cat}}
\newcommand{\bfCat}{\mathbf{Cat}}
\newcommand{\CAT}{\mathtt{CAT}}
\newcommand{\Cyl}{\mathtt{Cyl}}

\newcommand{\colim}{\operatorname{colim}}

\newcommand{\eq}{\mathit{eq}}

\newcommand{\Fin}{\mathit{Fin}}
\newcommand{\Fun}{\operatorname{Fun}}

\newcommand{\Map}{\operatorname{Map}}
\mathchardef\mhyphen="2D

\newcommand{\op}{\mathrm{op}}

\newcommand{\DOp}{\mathtt{DOp}}

\newcommand{\rlarrows}{\substack{\longrightarrow\\ \longleftarrow}}

\newcommand{\wt}{\widetilde}

\newcommand{\cC}{\mathcal{C}}
\newcommand{\cD}{\mathcal{D}}

\newcommand{\cL}{\mathcal{L}}

\newcommand{\cS}{\mathcal{S}}

\begin{document}

\title[]{Localization of lax symmetric monoidal categories}
\author{Vladimir Hinich}
\address{Department of Mathematics, University of Haifa, Mount Carmel, Haifa 3498838,  Israel}
\email{vhinich@gmail.com}
\author{Ieke Moerdijk}
\address{Department of Mathematics, University of Utrecht,
Budapestlaan 6, Utrecht,  Netherlands} 
 \email{I.Moerdijk@uu.nl}
 \begin{abstract}
Given a locally cocartesian fibration $f:X\to B$ and a collection of marked arrows $X^\circ\subset f^{-1}(B^\eq)$
in $X$ closed under locally cocartesian liftings, we prove that the localization $\cL(X,X^\circ)\to B$ is also a
locally cocartesian fibration whose fibers are localizations
of the fibers of $f$. This result is applied to the description of localizations of lax symmetric monoidal
categories.
\end{abstract}
\maketitle

\section{Introduction}

\subsection{}
In this note, we study localizations of a (co)lax symmetric monoidal $\infty$-category. Our motivation comes from the fact that one can obtain a symmetric monoidal structure on the category $\DOp$ of dendroidal $\infty$-operads by localizing a colax such structure on a category $P(\Phi)$ of presheaves, as described in Section 5.3 of our paper \cite{HM}. A lax symmetric monoidal category is defined as a Lurie $\infty$-operad that is simultaneously a locally cocartesian fibration over the category $\Fin_*$ of finite pointed sets. The underlying 
$\infty$-category appears in this description as the fiber at $\langle 1\rangle$ and we discuss the localization with respect to a collection
of marked arrows in this fiber.

One can think of two meaningful formulations of such localization. The first one is a ``fiberwise'' localization. We localize
the fiber at $\langle 1\rangle$ with respect to the marked arrows, and we are looking for a natural lax symmetric monoidal structure on it.
The second formulation requires that the localization of 
the total category of the fibration remains a locally cocartesian fibration that endows  the fiber at $\langle 1\rangle$
with a structure of lax symmetric monoidal category. 

The two formulations are a priori incomparable.  In general, if one just localizes the original category, it is perhaps not clear that it carries a lax symmetric monoidal category structure, i.e. that it can be extended to a locally cocartesian fibration. And if one localizes the total category of the fibration and shows it represents in fact a lax symmetric monoidal category, it is not clear that the underlying category of this localization is the intended localization of the original fibre over 
$\langle 1\rangle$. In our paper~\cite{HM} we implicitly used the first, fiberwise, localization, and we were not explicit about the relation
between the two.

 In this paper we show that,
for lax symmetric monoidal categories, assuming a natural compatibility
of the marking with the monoidal structure,  localization of the 
total category of the fibration provides a localization of the fibers. 
This implies in particular that, under the same compatibility condition, the total localization of a lax 
symmetric monoidal category provides the intended lax symmetric monoidal structure on 
the localization of the fiber at $\langle 1\rangle$. 
Another way of expressing this is by saying that taking fibers 
commutes with localization for suitably compatible markings in the fibres.
This property was already known for (non-lax) symmetric monoidal categories (see \cite{H.L}), but is easily seen to fail for operads in general.

\subsection{}

Recall that a marked category is a pair $(\cC,\cC^\circ)$ consisting
of an $\infty$-category $\cC$ and its subcategory $\cC^\circ$ 
containing the underlying space $\cC^\eq$. We denote by $\cC^\flat$
the minimally marked category $(\cC,\cC^\eq)$.
The localization $\cL(\cC,\cC^\circ)$ is a universal functor 
$\cC\to\cL(\cC,\cC^\circ)$ that carries the arrows of $\cC^\circ$ to equivalences.

In this note we are interested in localization of locally cocartesian fibrations. Let $f:\cC\to B$ be a locally cocartesian fibration.
This means that for any arrow
$\alpha:[1]\to B$ the base change $[1]\times_B\cC\to[1]$
is a cocartesian fibration. 
A {\sl marked locally cocartesian fibration} is a marked category 
$(\cC,\cC^\circ)$ endowed with a locally cocartesian fibration
$f:\cC\to B$ such that
\begin{itemize}
\item $f$ carries marked arrows in $\cC$ to equivalences in $B$.
\item for any $\alpha:b\to b'$ in $B$ the functor $\alpha_!:f^{-1}(b)\to f^{-1}(b')$ defined by the locally cocartesian lifting along $\alpha$,
preserves marked arrows.
\end{itemize}

The central result of this note is the following.
\begin{thm}
\label{thm:intro}(see Proposition~\ref{prp:loc})
Let $f:(\cC,\cC^\circ)\to B$ be a marked locally cocartesian fibration.
Then the localization $\cL(f):\cL(\cC,\cC^\circ)\to B$ is a locally 
cocartesian fibration whose fibers are the localizations of the fibers of 
$f$.
\end{thm}

\subsubsection{}
The proof of \ref{prp:loc} proceeds as follows.
First of all, we use straightening-unstraightening of locally cocartesian fibrations to deduce a weaker form of our localization theorem, 
see~\ref{prp:lccfloc}. This weaker form claims that there is a universal
map of locally cocartesian fibrations $\cC\to\wt\cC$ over $B$ that carries
the arrows of $\cC^\circ$ to equivalences and induces a fiberwise localization, but leaves open the relation to the localization of the total category.
This weaker claim is sufficient to justifying our claim 
in Proposition 5.3.1 of \cite{HM}.

Next, we deduce Proposition~\ref{prp:loc} for the special case $B=[n]$,
by induction in $n$, using an inductive description of locally cocartesian fibrations over $[n]$. Finally, we use a standard presentation of $B$ as a colimit in order to deduce~\ref{prp:loc} for general $B$, and hence 
prove Theorem~\ref{thm:intro}.

\subsection{}
\label{ss:laxSM}
A lax symmetric monoidal (SM) category can be defined as an operad $p:\cC\to\Fin_*$ that is,
in addition,
a locally cocartesian fibration. Recall that an operad is called a symmetric monoidal category
if it is a cocartesian fibration. In particular, a lax SM category is 
symmetric monoidal if and only if the composition of locally cocartesian liftings is
locally cocartesian.
We denote by $\cC_1$ the fiber of $p$ at $\langle 1\rangle\in\Fin_*$.
This is the category underlying the lax SM category $\cC$.

Applying \ref{prp:lccfloc} and \ref{prp:loc} to lax symmetric monoidal categories, we deduce Theorem~\ref{thm:laxSMloc} below.

Let $p:\cC\to\Fin_*$ be a lax SM category structure on $\cC_1=\{\langle 1\rangle\}\times_{\Fin_*}\cC$. Given a marking $\cC_1^\circ$ on  $\cC_1$,
we define the marking $\cC^\circ$ on $\cC$ {\sl generated by }
$\cC_1^\circ$ by the formula $\cC^\circ=\sqcup\cC_n^\circ$ where
$\cC_n^\circ\subset\cC_n=\{\langle n\rangle\}\times_{\Fin_*}\cC$ is the 
marking  that makes the standard equivalence $\cC_n\to\cC_1^n$ an 
equivalence of marked categories.

\begin{thm}
\label{thm:laxSMloc}
Let $\cC$ be a lax SM category. Let $\cC^\circ_1\subset\cC_1$ be a subcategory containing $\cC_1^\eq$ such that the 
marking $\cC^\circ$ of $\cC$ generated by $\cC^\circ_1$ makes $(\cC,\cC^\circ)$ into a marked locally cocartesian fibration over $\Fin_*$
as defined above.
Then the localization $\cL(\cC,\cC^\circ)\to\Fin_*$ is a lax SM category structure on  $\cL(\cC_1,\cC^\circ_1)$. 
\end{thm}

\subsection{Contents of the sections} In Section~\ref{sec:strunstr} we review the straightening and unstraightening for locally cocartesian fibrations.  In Section 3 we prove, using straigtening--unstraightening, 
Proposition~\ref{prp:lccfloc} that presents a fiberwise localization
of a locally cocartesian fibration. Finally, in Section~\ref{sec:strong}
we present a stronger version claiming that the universal localization described in~\ref{prp:lccfloc} is given by a  localization of the total category of fibration.

\subsubsection{Acknowledgements}We are grateful to 
 Manuel Krannich and Sander Kupers for pointing out to us the ambiguity 
 of our treatment of localization in~\cite{HM}, 5.3.

\section{Straightening-unstraightening for locally cocartesian fibrations}
\label{sec:strunstr}
The weak version of the localization result,
Proposition~\ref{prp:lccfloc}, is based on  a  version of
straightening/unstraightening equivalence for locally cocartesian 
fibrations that  we present in Proposition~\ref{prp:lccf}. This result
is based on the study of double categories and $2$-categories
that we   recall in Section~\ref{ss:double}. In particular, we  define 
the $2$-categories $\bfCat$ and $\bfCat^+$ of $\infty$-categories and of marked $\infty$-categories used in the sequel.  

\subsection{Double categories, $2$-categories}
\label{ss:double}

\subsubsection{}

A double category is, by definition, a Segal object in $\Cat$, that is,
a simplicial object $\Delta^\op\to\Cat$ satisfying
the Segal condition. We will present double categories by the
corresponding cartesian fibrations $p:X\to\Delta$. Functors between double 
categories are defined as maps of cartesian fibrations preserving
the cartesian arrows.

An oplax functor $\phi:X\to Y$ between two double categories is a 
functor $\phi:X\to Y$ over $\Delta$ preserving cartesian liftings of the 
inerts (that is, injective maps $a:[k]\to [n]$ given by the formula
$a(i)=r+i$ for fixed $r$). An oplax functor is called unital if, in addition, it preserves the cartesian liftings of the arrows
$[n]\to[0]$. (Note that by the Segal condition the same will be true for any surjective map
$[n]\to[k]$.)

\subsubsection{Variants}

A double category $X_\bullet$ is called  a 
{\sl $2$-precategory} if $X_0$ is a space.  
	A {\sl $2$-category} is a $2$-precategory additionally satisfying the completeness condition.   
Given a double category $X_\bullet$, the induced   $2$-precategory is obtained by replacing $X_0$ with $X_0^\eq$ and $X_n$ with the preimage
of $(X^\eq_0)^{n+1}$ under the map $X_n\to X_0^{n+1}$.

\subsubsection{A category  as a $2$-category}
\label{sss:cat-2cat}

Any category $\cC$ gives rise to a Segal space $\Delta^\op\to\cS$, so,
it defines a Segal object in categories $\Delta^\op\to\cS\to\Cat$.
In other words, any category $\cC$  defines a  double category that is, in fact, a $2$-category. If $\cC$ is a category, the cartesian presentation of the corresponding $2$-category is $p:\Delta_{/\cC}\to\Delta$, where 
$\Delta_{/\cC}$ is the usual category of simplices in $\cC$.

\subsubsection{The $2$-category $\bfCat$ of categories}
\label{sss:cat}
We define the double category of categories as a cartesian fibration
$p:\CAT\to\Delta$ as follows. Its total category $\CAT$ is the full subcategory of $\Fun([1],\Cat)$ spanned by the cartesian fibrations
$\alpha:X\to [n]^\op$
(of course, $[n]^\op$ is canonically isomorphic to $[n]$). The projection
$p:\CAT\to\Delta$ carries $\alpha$ to its codomain. The functor $p$
is a cartesian fibration as $\Cat$ has fiber products and cartesian fibrations are closed under base change.  The cartesian fibration
$p:\CAT\to\Delta$ satisfies the Segal condition, so this is a double category.

We denote by $\bfCat$ the $2$-precategory induced by $\CAT$. It is a 
$2$-category.

\subsubsection{The $2$-category $\bfCat^+$ of marked categories}

Marked categories, that is, pairs $(\cC,\cC^\circ)$ consisting of a category $\cC$ and a subcategory $\cC^\circ\supset\cC^\eq$, and marked functors, form a category denoted by $\Cat^+$.
There is also a $2$-categorical version that we will now describe.
We define $\CAT^+$ as the full subcategory of $\Fun([1],\Cat^+)$ spanned by the cartesian fibrations
$\alpha:X\to  [n]^{\op\flat}$ (the category $[n]^\op$ with no nontrivial markings)  for which the induced functors $\alpha^{-1}(j)\to \alpha^{-1}(i)$ preserve the markings. The functor $p:\CAT^+\to\Delta$
is a cartesian fibration and the Segal condition obviously holds.
The $2$-precategory $\bfCat^+$ is the one induced by 
$\CAT^+$. It is a $2$-category.  

There is an obvious functor $\bfCat^+\to\bfCat$ over $\Delta$ forgetting the marking. It preserves the cartesian arrows, so it defines a functor between the $2$-categories.

\subsection{}
The straightening--unstraightening result for 
 locally cocartesian fibrations referred to in 1.2.2 is the following.

\begin{prp}
\label{prp:lccf}
There is a natural equivalence between the category of locally cocartesian 
fibrations over $B\in\Cat$ and the category of oplax unital functors $B\to\bfCat$ to the $2$-category of categories. 
\end{prp}

A version of Proposition~\ref{prp:lccf} in the framework of scaled simplicial sets is proven in \cite{LG}, Theorem 3.8.1. Another approach is
based on the straightening-unstraightening equivalence for
cartesian fibrations over $n$-categories established by J.~Nuiten, see~\cite{N}, Theorem 6.1.
See also~\cite{AMGR}, Theorem B.3.7 and B.1.23.

\section{Localization of locally cocartesian fibrations, I}

\subsection{Marked locally cocartesian fibrations}
In this subsection we interpret marked locally cocartesian fibrations
over $B$ as oplax unital functors from $B$ to $\bfCat^+$.

A marked locally cocartesian fibration is a map $p:(\cC,\cC^\circ)\to B^\flat$
of marked categories such that the underlying locally cocartesian
fibration $p:\cC\to B$ has the property that
for any $\phi:a\to b$ in $B$ the functor $\phi_!:p^{-1}(a)\to p^{-1}(b)$ 
preserves the markings. 
In this subsection we show that a marked locally cocartesian fibration
$p:(\cC,\cC^\circ)\to B^\flat$
defines canonically an oplax unital functor 
$B\to\bfCat^+$ such that the composition with the forgetful functor 
$B\to\bfCat^+\to\bfCat$ classifies the locally cocartesian fibration
$p:\cC\to B$.

Let $\cC\in\Cat$.   Given $f:\cC\to\bfCat$, we want to describe possible
liftings of $f$ to a functor $\cC\to\bfCat^+$.
Recall the following easy fact.
\begin{lem} 
\label{lem:sub}
Let $f:\cC\to\cD$ be a functor. Let for each $x\in\cC$
a subobject $e_x$ of $f(x)\in\cD$ is given, so that for any 
$\alpha:x\to y$ in $\cC$ the composition $e_x\to f(x)\to f(y)$
factors through $e_y$. Then the collection of objects $e_x$ defines
a unique subfunctor $g:\cC\to \cD$ with canonical equivalences
$g(x)=e_x$.
\end{lem}
\begin{proof}
By composing $f$ with the Yoneda embedding $\cD\to P(\cD)$
and transposing to a functor $\cC\times\cD^\op\to\cS$, we may reduce the
problem to the case $\cD=\cS$, where the result follows easily by 
unstraightening.
\end{proof}

\begin{prp}
\label{prp:lifting}
Given a functor $f:\cC\to\bfCat$ and a collection of subcategories
$W_x\subset f(x)\in\Cat$ for each $x\in\cC$ containing $f(x)^\eq$, such 
that for any $\phi:x\to y$ one has $f(\phi)(W_x)\subset W_y$, there 
exists a unique lifting of $f$ to a functor $\hat f:\cC\to\bfCat^+$
carrying $x\in\cC$ to a marked category $(f(x),W_x)$.
\end{prp}
\begin{proof}
Since $\bfCat^+$ and $\bfCat$
are subcategories of $\Fun([1],\Cat^+)$ and of $\Fun([1],\Cat)$
respectively, the functor $f:\cC\to\bfCat$ gives rise to a functor
$f':\cC\times[1]\to\Cat$. The collection of $W_x$ allows one to
define a subfunctor $\cC\times[1]\to\Cat$ carrying $(x,0)$ to $W_x\subset
f(x)$ and $(x,1)$ to $p\circ f(x)$. By~\ref{lem:sub} this is equivalent to a functor 
$\cC\times[1]\to\Cat^+$ lifting $f$. 
\end{proof}

\subsubsection{}
Following~\ref{sss:cat-2cat},
an oplax unital functor from $B$ to $\bfCat$ is a functor
$$
F:\Delta_{/B}\to\bfCat
$$
over $\Delta$  preserving cartesian lifts of inerts and of $[n]\to[0]$.
 
Let $p:(\cC,\cC^\circ)\to B^\flat$ be a marked locally cocartesian fibration\and let $F:\Delta_{/B}\to\bfCat$ be the unital oplax functor classifying (the unmarked version of) $p$. 
Let  $s:[n]\to B$ in $\Delta_{/B}$ be an $n$-simplex
$$
s: a_0\stackrel{\phi_1}{\to}\ldots\stackrel{\phi_n}{\to}a_n.
$$
The category $F(s)\in\bfCat_n$
is the cartesian fibration over $[n]^\op$ classified by a sequence of
functors denoted
$$
\cC_{a_0}\stackrel{(\phi_1)_!}{\to}\ldots\stackrel{(\phi_n)_!}{\to}\cC_{a_n}.
$$
The collection of objects of $F(s)$ is the union of all objects of 
the categories $\cC_{a_i}$. We define $W(s)$ as the subcategory of $F(s)$ spanned by the arrows that belong to one of $\cC^\circ_{a_i}$. The collection of $W(s)$ satisfies
Proposition~\ref{prp:lifting}, so it allows to canonically factor
$F:\Delta_{/B}\to\bfCat$ as a composition
 $\Delta_{/B}\stackrel{F^+}{\to}\bfCat^+\to\bfCat$.
 The functor $F^+:\Delta_{/B}\to\bfCat^+$ is the oplax unital functor
 classifying the marked lax cocartesian fibration $p:(\cC,\cC^\circ)\to B^\flat$.

\subsection{Localization}
\subsubsection{}

The localization functor $L:\Cat^+\to\Cat$ extends to a functor 
$$
\Fun([1],\Cat^+)\to\Fun([1],\Cat).
$$

We denote by $L$ its restriction to the full subcategory
$\bfCat^+$ of  $\Fun([1],\Cat^+)$.
By~\cite{H.L}, 2.1.4, the localization functor commutes with the projections to $\Delta$ and preserves  cartesian fibrations, so it  carries $\bfCat^+\subset\Fun([1],\Cat^+)$ to $\bfCat\subset
\Fun([1],\Cat)$ and
defines a functor between these $2$-categories.

The functor $L$ is left adjoint to $\flat:\bfCat\to\bfCat^+$
sending a category $B$ to the pair $B^\flat=(B,B^\eq)$; this 
reflects the universality of localization.

\subsection{Localization of lccf, a weaker form}
We can now easily deduce a weak form of our localization result for
locally cocartesian fibrations.
 
\begin{prp}
\label{prp:lccfloc}
Let $p:(\cC,\cC^\circ)\to B^\flat$ be a marked locally cocartesian fibration, such that for any $\alpha:x\to y$ in $B$ the functor
$\alpha_!:p^{-1}(x)\to p^{-1}(y)$ preserves the markings. Then 
there is a universal morphism of locally cocartesian fibrations
$$
\cC\to\cL
$$
over $B$
that carries all the arrows in $\cC^\circ$ to equivalences. Moreover,
for each $x\in B$ the induced map $\cL(\cC_x,\cC^\circ_x)\to\cL_x$ is an equivalence.
\end{prp}
\begin{proof}
The locally cocartesian fibration $p:\cC\to B$ gives rise to an oplax unital functor $B\to \bfCat$. The marking on $\cC$ gives rise to its canonical lifting to an oplax unital functor $B\to\bfCat^+$. Composing
it with the localization functor $\cL:\bfCat^+\to\bfCat$, we get an oplax 
unital functor $B\to\bfCat$ that defines a locally cocartesian fibration  
$\cL\to B$. By construction, the fiber $\cL_x$, $x\in B$, 
is the localization of $\cC_x$ with respect to $\cC_x^\circ$.
The universal property of $\cC\to\cL$ 
follows from the fact that the functor $L:\bfCat^+\to\bfCat$ is left adjoint to $\flat:\bfCat\to\bfCat^+$.
 
\end{proof}

\section{Localization of locally cocartesian fibrations, II}
\label{sec:strong}

\subsection{Locally cocartesian fibrations over $[n]$}

In this subsection we provide a recursive description of locally cocartesian fibrations over $[n]$. This will allow us to deduce
Proposition~\ref{prp:loc} that strengthens \ref{prp:lccfloc}, in the special case $B=[n]$. This will eventually lead to a proof of
Proposition~\ref{prp:loc} in general.

\subsubsection{}
Recall that, for $f:X\to Y$ in $\Cat$, the cylinder of $f$ is defined
as the pushout
$$
\Cyl(f)=(X\times[1])\sqcup^XY,
$$
where the map $X\to X\times[1]$ is induced by $\{1\}\to[1]$.
 It is well-known (but not at all obvious) that $\Cyl(f)\to[1]$ is the cocartesian fibration classified by the functor $[1]\to\Cat$ defined by
$f:X\to Y$. 

\subsubsection{}

Let $p:X\to[n]$ be a category over $[n]$. We denote by
$p':X'\to[n-1]$ the base change of $p$ with respect to
the map $\partial^0:[n-1]\to[n]$. We denote $X_0:=p^{-1}(0)$.
\begin{lem} The map $p$ is a locally cartesian fibration iff
\begin{itemize}
\item[1.] $p'$ is a locally cartesian fibration.
\item[2.] The map $X'\to P(X_0)$ obtained by restriction the Yoneda embedding $X\to P(X)$, factors through the Yoneda embedding $X_0\to P(X_0)$.
\end{itemize}
\end{lem}
\begin{proof}Condition 2 means that for any $x'\in X'$ its image in $P(X_0)$ is representable. This immediately implies the claim.
\end{proof}

We will apply the above lemma to locally cocartesian fibrations.
Let $p:X\to[n]$ be a category over $[n]$ and let $p':X'\to[n-1]$
be obtained from $p$ by the base change along $\partial^n:[n-1]\to[n]$.
Let, furthermore, $q:[n]\to[1]$ be defined by the formula
$q(n)=1$, $q(i)=0$ for $i<n$. 

\begin{crl}The map $p$ is locally cocartesian fibration iff
\begin{itemize}
\item[1.] $p'$ is locally cocartesian fibration.
\item[2.] $q\circ p:X\to[1]$ is a cocartesian fibration
(classified by a map $(q\circ p)^{-1}(0)=X'\to X_n$).
\end{itemize}
\end{crl}
\begin{crl}
\label{crl:cyl-expression}
Let $p:X\to[n]$ be a locally cocartesian fibration. Then
one has a natural equivalence $X=\Cyl(\phi:X'\to X_n)$ where $\phi$
is the map classifying the cocartesian fibration $q\circ p$.
\end{crl}

Corollary~\ref{crl:cyl-expression} allows one to describe a locally cocartesian fibration over $[n]$ as an iterated cylinder.

\subsubsection{Example, $n=[2]$}
Let $f_{ij}:X_i\to X_j$ classify the cocartesian fibration obtained by
the base change along $s:[1]\to[2]$, $s(0)=i$ and $s(1)=j$.
We have $X'=\Cyl(f_{01})$ so that a map $\phi:X'\to X_2$ is given by
a pair of maps $f_{02}:X_0\to X_2$, $f_{12}:X_1\to X_2$ and a homotopy
from $f_{02}$ to $f_{12}\circ f_{01}$ given by $\phi:X_0\times[1]\to X_2$.
We finally have $X=\Cyl(\phi:\Cyl(f_{01})\to X_2)$.

\subsection{A colimit presentation}
 
 In this subsection we  present a category $p:X\to B$ over $B$ as a 
colimit over the simplices of $B$.

Let $\Delta_{/B}$ be the category of simplices in $B$. The morphism
$p:X\to B$ gives rise to a functor $\phi_X:\Delta_{/B}\to\Cat_{/B}$
carrying $s:[n]\to B$ to the fiber product $\phi_X(s)=[n]\times_BX$.

The functor $\phi_X$ canonically extends to 
$\tilde\phi_X:(\Delta_{/B})^\triangleright\to\Cat_{/B}$ carrying the terminal object of $(\Delta_{/B})^\triangleright$ to $p:X\to B$.

\begin{lem}
\label{lem:colim}
The functor $\tilde\phi_X$ presents $p:X\to B$ as the
colimit of $\phi_X$.
\end{lem}
\begin{proof} We start with an elementary observation. Let $\cC$ be a category and let $Y:\cC\to P(\cC)$ denote the Yoneda embedding. We complete $Y$ to a functor $\tilde Y:\cC^\triangleright\to P(\cC)$ carrying the terminal object
of $\cC^\triangleright$ to a terminal object $t\in P(\cC)$.
$\tilde Y$ is well-known to be a colimit diagram. Since
the functor $\times F$, $F\in P(\cC)$, has right adjoint, 
the functor $\tilde Y\times F$ is also a colimit diagram.

We are almost done. The category $\Cat_{/B}$ can be presented 
as a Bousfield localization of $P(\Delta_{/B})$ as follows.

The full embedding
$R:\Cat\to P(\Delta)$ has a left adjoint localization functor $L:P(\Delta)\to\Cat$.
Given $B\in\Cat$, we get an adjuction
$$
L:P(\Delta)_{/R(B)}\rlarrows \Cat_{/B}:R.
$$
 The category $P(\Delta)_{/R(B)}$  is naturally equivalent to the category of presheaves $P(\Delta_{/B})$ where 
 $\Delta_{/B}$ is the category of simplices in $B$; see, for instance,~\cite{HM}, 2.3.3. We will prove that the composition $R\circ\tilde\phi_X$ is a colimit 
diagram. This implies the claim of the lemma as the localization $L$ preserves colimits and
$\tilde\phi=L\circ R\circ\tilde\phi$.

Let $F=R(X)$. We have $R\circ\tilde\phi_X(s)=Y(s)\times F$,
so by the elementary observation at the beginning of the proof, $R\circ\tilde\phi_X$ is a colimit diagram.

\end{proof}

\subsection{Localization of lccf, a stronger form}

Proposition~\ref{prp:loc} below is the central result of this note.
It is used in~\ref{ss:loclaxSM} in the proof of
Theorem~\ref{thm:laxSMloc}.

\begin{prp}
\label{prp:loc}
Let $p:(X,X^\circ)\to B^\flat$ be a marked locally cocartesian fibration. This means that
\begin{itemize}
\item[1.] $p:X\to B$ is a locally cocartesian fibration.
\item[2.] for any arrow $f$ in $X^\circ$ its image  $p(f)$ is an equivalence in $B$.
\item[3.] For any arrow $\alpha$ in $B$ the functor $\alpha_!$ preserves
marked arrows in $X$.
\end{itemize}
Then the localization $\cL(p):\cL(X,X^\circ)\to B$ is a locally cocartesian fibration whose fibers are the localizations of the fibers of $p$.
\end{prp}

We will first prove the result in the case $B=[n]$ and
then deduce the general case using the colimit presentation of Lemma~\ref{lem:colim}.
\subsubsection{}
\label{sss:1}
If $B=[1]$, $p:X\to B$ is a cocartesian fibration. If $X_i=p^{-1}(i)$, 
let $\phi:X_0\to X_1$ classify $p:X\to[1]$. Let $X^\circ_i=X^\circ\cap X_i$, $i=0,1$. 
The natural map $\Cyl(\phi)\to X$ is an equivalence. Any map $q:X\to Z$
carrying $X^\circ$ to $Z^\eq$, uniquely factors through $\Cyl(\psi)$
where $\psi:\cL(X_0,X^\circ_0)\to\cL(X_1,X^\circ_1)$ is obtained from 
$\phi$ by localization. This implies that $\cL(X,X^\circ)=\Cyl(\psi)$.
This proves the claim for $B=[1]$.

\subsubsection{}
\label{sss:simplex}
We   prove the  proposition for $B=[n]$ by induction
in $n$. Given $p:X\to[n]$ and $p':X'\to[n-1]$ as 
in~\ref{crl:cyl-expression},
we have $X=\Cyl(\phi:X'\to X_n)$. If $p:(X,X^\circ)\to[n]^\flat$ is a marked locally cocartesian fibration, $(X',X'\cap X^\circ)\to[n-1]^\flat$
is a marked locally cocartesian fibration by the inductive hypothesis, 
so by~\ref{sss:1} the claim holds for $(X,X^\circ)$.

\subsubsection{$B$ is general}

We first of all apply Proposition~\ref{prp:lccfloc} to get 
a fiberwise localization of $X$, that is a morphism of lccf
$X\to\wt X$ over $B$ that carries $X^\circ$ to   $(\wt X)^\eq$ and such that the canonical maps $\cL(X_b, W\cap X_b)\to \wt X_b$ are equivalences for all $b\in B$. We have to prove that the canonical map 
$L:\cL(X,X^\circ)\to\wt X$ is an equivalence.  Note that the construction of
the fiberwise localization presented in the proof of Proposition~\ref{prp:lccfloc}, commutes with the base change, so
that for any $s:[n]\to B$ the base change $\wt X_s=[n]\times_B\wt X$ is the fiberwise
localization of $p_s:X_s\to[n]$ with respect to $X^\circ_s=[n]\times_BX^\circ$. By~\ref{sss:simplex}, the natural map
$$
\cL(X_s,X_s^\circ)\to \wt X_s
$$
is an equivalence.

We will use the colimit presentation 
of Lemma~\ref{lem:colim} for $X$ and for $\wt X$.  We have
a commutative diagram

\begin{equation}
\xymatrix{
& \colim_{s\in\Delta_{/B}}X_s\ar[d] \ar[r] & X \ar[d]\\
&\colim_{s\in\Delta_{/B}}\cL(X_s,W_s)\ar^\alpha[d]\ar[r] &\cL(X,W) \ar^\beta[d] \\
&\colim_{s\in\Delta_{/B}}\wt X_s\ar[r]& \wt X 
},
\end{equation}
where the top and the bottom horizontal arrows are equivalences
by Lemma~\ref{lem:colim}. The middle horizontal arrow is also equivalence
by the universality of localization. Since $\alpha$ is an equivalence, 
$\beta$ is also an equivalence.
 
This proves Proposition~\ref{prp:loc}.

\subsection{Localization of lax symmetric monoidal categories}
\label{ss:loclaxSM}

The claim about localization of lax SM categories is an easy  consequence of Proposition~\ref{prp:loc}.

\begin{lem}
\label{lem:lcc-cc}
Let $p:X\to B$ be a locally cocartesian fibration and $\alpha:a\to b$
be an arrow in $B$. The following conditions are equivalent.
\begin{itemize}
\item[1.] Any locally cocartesian lifting of $\alpha$ is (globally) cocartesian.
\item[2.] For any commutative triangle $\gamma=\beta\circ\alpha$ in $B$
the natural map $\gamma_!\to\beta_!\circ\alpha_!$ is an equivalence.
\end{itemize}
\end{lem}
\begin{proof}This is a variation of Lemma 2.4.2.7 of~\cite{L.T}.
The implication $(1)\Longrightarrow (2)$ immediately follows from
{\sl loc. cit.} Assume that $(2)$ holds and let $\bar\alpha$ be a 
lcc lifting of $\alpha$. Let $\bar\gamma=\bar\beta\circ\bar\alpha$ be lcc.
We have to show that $\bar\beta$ is also lcc. In fact, choose $\bar\beta'$
over $\beta$ that is lcc. One has a decomposition 
$\bar\beta=u\circ\bar\beta'$ where $p(u)$ is an equivalence. A composition
$\bar\gamma':=\bar\beta'\circ\bar\alpha$ is a lcc lifting of $\gamma$,
as well as $\bar\gamma=\bar\gamma'\circ u$. This implies that $u$ is
an equivalence.
\end{proof}

\subsubsection{Proof of Theorem~\ref{thm:laxSMloc}}

Let $p:\cC\to\Fin_*$ be a lax SM category. Let $\cC^\circ_1\subset\cC_1$
be a marking closed under the multiple tensor products and let $\cC^\circ$ be the marking of $\cC$ generated by $\cC_1^\circ$ as defined 
in~\ref{ss:laxSM}.
 
Let $\cL=\cL(\cC,\cC^\circ)$. By Proposition~\ref{prp:loc}, $\cL\to\Fin_*$
is a locally cocartesian fibration.
It remains to verify that the map $\cL\to\Fin_*$ defines an operad.
This amounts to the following:
\begin{itemize}
\item[1.] Inerts in $\Fin_*$ have cocartesian lifting.
\item[2.] The cocartesian liftings of the inerts
$\rho^i:\langle n\rangle\to\langle 1\rangle$
induce equivalences $\cL_n\to\cL_1^n$.
\item[3.] For any $x,y\in\cL$ with $y\in\cL_n$ and $y=\oplus y_i$
and $f:p(x)\to p(y)=\langle n\rangle$ in $\Fin_*$, 
the map  
\begin{equation}
\label{eq:decompose}
\Map^f(x,y)\to\prod_i\Map^{\rho^i\circ f}(x,y_i)
\end{equation}
is an equivalence.
\end{itemize}

Let us verify these properties.

1. By Lemma~\ref{lem:lcc-cc}, the morphism of locally cocartesian fibrations
$\cC\to\cL$ preserves the property of an arrow $\alpha$ to have cocartesian liftings. Therefore, since inerts in $\Fin_*$ admit cocartesian liftings in $\cC$, the same holds for $\cL$.

2. Choose cocartesian liftings for $\rho^i:\langle n\rangle\to\langle 1\rangle$.  We have an equivalence $\cC_n\to\cC_1^n$ that preserves the marked arrows. Since the localization preserves products, we get an equivalence $\cL_n\to\cL_1^n$.

3. In the case when $\cC$ is a lax SM category, the condition
(\ref{eq:decompose}) is equivalent to the condition 
that the map $(\rho^i\circ f)_!\to\rho^i_!\circ f_!$ is an equivalence.
Since the functor $\cC\to\cL$ preserves locally cocartesian arrows,
the same condition holds for $\cL$.

\end{document}